%
\magnification=\magstep1
\input amstex
\UseAMSsymbols
\input pictex
\NoBlackBoxes
\font\kl=cmr8  \font\gross=cmbx10 scaled\magstep1 
   \font\rmk=cmr8    \font\itk=cmti8  \font\ttk=cmtt8

   \newcount\notenumber
   
   \def\note{\advance\notenumber by 1 
       \plainfootnote{$^{\the\notenumber}$}}  

\def\Mminus{M(n\!-\!1)}

\def\Hom{\operatorname{Hom}}

\def\rad{\operatorname{rad}}
\def\soc{\operatorname{soc}}

   

\centerline{\gross Indecomposables  live in all smaller lengths.}     
		   		       \medskip\smallskip
\centerline{Claus Michael Ringel}
		  	    \medskip\smallskip
\plainfootnote{}
{\rmk 2000 \itk Mathematics Subject Classification. \rmk 
Primary 
        16D90, 
        16G60. 
Secondary:
        16G20. 
}

{\narrower \kl Abstract.  Let $\ssize \Lambda$ be a 
finite-dimensional $\ssize k$-algebra
with $\ssize k$ algebraically closed. Bongartz has recently shown that the existence of 
an indecomposable $\ssize \Lambda$-module 
of length $\ssize n > 1$ implies that also indecomposable $\ssize \Lambda$-modules of length 
$\ssize n-1$ exist. Using a slight modification of his arguments, we strengthen the 
assertion as follows: If there is an indecomposable module of length $\ssize n$, then there is also 
an accessible one. Here, the accessible modules are defined inductively, as follows:
First, the simple modules are accessible.
Second, a module of length $\ssize n \ge 2$ is accessible provided it is indecomposable and 
there is a submodule or a factor module of length $\ssize n-1$ which is accessible. 
\par}

	\bigskip

Let $k$ be an algebraically closed field. Let $\Lambda$ be a finite-dimensional $k$-algebra,
we may (and will) assume that $\Lambda$ is basic.
We are interested in 
(usually finite-dimensional left) $\Lambda$-modules. A recent preprint [B3]
of Bongartz with the same title is devoted to a proof of the following important result:
       \medskip
{\bf Theorem (Bongartz 2009).} {\it Let $\Lambda$ be a finite-dimensional $k$-algebra
with $k$ algebraically closed. If there exists an indecomposable $\Lambda$-module 
of length $n > 1$, then there exists an indecomposable $\Lambda$-module of length $n-1$.}
   \medskip
Unfortunately, the statement does not assert any relationship between 
the modules of length $n$ and 
those of length $n-1$. There is the following open problem:
{\it Given an indecomposable $\Lambda$-module $M$ of length $n\ge 2$. 
Is there an indecomposable submodule or factor module of length $n-1?$}
	
{\bf Remarks. (1)} This is the case for $\Lambda$ being representation-finite or 
tame concealed, as Bongartz [B1, B2] has shown 
already in 1984 and 1996, respectively, but the answer is unknown in general.
A positive answer would have to be
considered as a strong finiteness condition --- after all, if we consider 
for example any quiver of type $\Bbb A_\infty^\infty$, then there is a unique
minimal faithful representation $M$, it is indecomposable, but all its
maximal submodules as well as all the 
factor modules $M/S$ with $S$ simple, are decomposable. 
	
{\bf (2)} It is definitely necessary to look both for submodules
and factor modules, since for suitable algebras $\Lambda$, 
there are indecomposable modules $M$ which are not simple and have
no maximal submodules which are indecomposable. 
Any local module of length at least 3 and
Loewy length 2 is an example. And dually, there are indecomposable 
modules $M$ of length $n\ge 3$ such that all factor modules of length $n-1$ are decomposable.
	       
{\bf (3)} In case we weaken the assumption on the base field $k$, then we may find
counter-examples. For instance, let $k$ be the field with 2 elements, 
$Q$ the 3-subspace quiver (this is the quiver of type $\Bbb D_4$ with one sink and
3 sources) and $M$ the (unique) indecomposable $kQ$-module of length 5. There is also only
one indecomposable $kQ$-module $N$ of length 4. Now $N$ cannot be a submodule of $M$, since
we even have  $\Hom(N,M) = 0$. But $N$ is also not a factor module of $M$, since
$\Hom(M,N)$ is a 2-dimensional $k$-space and the three non-zero elements in $\Hom(M,N)$
all have images of length 3. For dealing with an arbitrary field $k$, one may
ask: {\it Given an indecomposable $\Lambda$-module $M$ of length $n\ge 2$,
is there an indecomposable module $N$ of length $n-1,$ generated or 
cogenerated by $M$ ?}
	\medskip
The present note modifies slightly the arguments of Bongartz in [B3] in order to
strengthen his assertion. We define inductively 
{\it accessible} modules: First, the simple modules are accessible.
Second, a module of length $n \ge 2$ is accessible provided it is indecomposable and 
there is a submodule or a factor module of length $n-1$ which is accessible. 
The open problem mentioned above can be reformulated as follows: Are all indecomposable modules
accessible? For a certain class of algebras, we are going to construct 
a suitable number of accessible modules of arbitrarily large length.
	    
We call an inclusion of modules $M' \subseteq M$ {\it uniform}, provided any submodule
$U$ with $M' \subseteq U \subseteq M$ is indecomposable (this is related to the 
well-accepted notion of a uniform module: a module $M$ is uniform provided it is non-zero
and any inclusion $M'\subset M$ with $M'\neq 0$ is uniform). If $M' \subseteq M$ is a uniform
inclusion, then $\soc M' = \soc M$. The converse is not true: the inclusion of a module $M'$
into its injective envelope $E(M')$ is uniform only in case $M'$ itself is uniform, however $M'$
and $E(M')$ always have the same socle. 
If $M' \subseteq M$ is a uniform inclusion, then the module $M$ is obtained from the
indecomposable module $M'$ by successive extensions (from above) using simple modules, with
all the intermediate modules being indecomposable. In particular, if 
$M' \subseteq M$ is a uniform inclusion and $M'$ is accessible, then also $M$
is accessible. There is the dual notion of a
couniform projection: If $X$ is a submodule of $M$, then the canonical map $M \to M/X$
is said to be a {\it couniform projection} provided all the modules $M/X'$ with $X'$ a submodule
of $X$ are indecomposable. Of course, if $M \to M''$ is a couniform projection and
$M''$ is accessible, then also $M$ is accessible.
   
Our aim is to show that all representation-infinite algebras have accessible modules
of arbitrarily large length. As Bongartz has pointed out (see 
the proof of the Corollary below), it is actually enough to look at 
non-distributive algebras. We recall that a finite-dimensional algebra is said to be
{\it non-distributive} in case its ideal lattice is not distributive. 
     \bigskip
{\bf Theorem.} {\it Let $\Lambda$ be a non-distributive algebra.
Then there are $\Lambda$-modules $M(n),$ $R(n),$ $W(n)$ and non-invertible 
homomorphisms
$$
\align
 W(1) &\leftarrow R(2) \ \leftarrow \ M(2)\ \to\quad R(3)\quad \to \quad  W(3)\quad  \leftarrow \dots \cr 
\dots \to 
 W(2n-1) &\leftarrow R(2n) \leftarrow M(2n) \to R(2n\!+\!1) \to W(2n\!+\!1) \leftarrow \dots \cr 
\endalign
$$
where the arrows pointing to the left are couniform projections and those pointing
to the right are uniform inclusions, and such that $W(1)$ is a uniform module.}
   \medskip 
By induction it follows that all these modules $M(n), R(n), W(n)$ are accessible. 
In particular, we see that {\it a non-distributive algebra $\Lambda$
has accessible modules of arbitrarily large
length.}
	\medskip
It seems to be surprising that here we deal with a very natural 
question that had not yet been settled for non-distributive algebras. 
Note that the 
class of non-distributive algebras
was the first major class of representation-infinite algebras studied in
representation theory, see Jans [J], 1957. 
Before we turn to the proof of the Theorem, let us derive the following consequence.
       \bigskip
{\bf Corollary.} {\it Let $\Lambda$ be a finite-dimensional $k$-algebra
with $k$ algebraically closed.
If there is an indecomposable module of length $n$, then there is  
an accessible one of length $n$.}
    \medskip
Proof of Corollary. As we have mentioned, for a representation-finite algebra all the indecomposable
modules are accessible, thus we can assume that $\Lambda$ is representation-infinite.
According to Roiter's solution [R] of the first Brauer-Thrall conjecture, a
representation-infinite algebra has indecomposable modules of arbitrarily large length,
thus we have to show that $\Lambda$ has accessible modules of any length.
Clearly, we can assume that $\Lambda$ is minimal representation-infinite (this means that $\Lambda$ is
representation-infinite and that any proper factor algebra is representation-finite).
			
According to Bongartz [B3, section 3.2] we only have to consider algebras with
non-distributive ideal lattice: Namely, if $\Lambda$ is minimal representation-infinite
and the ideal lattice of $\Lambda$ is distributive, then
the universal cover is interval-finite and the fundamental group
is free; using covering theory, the problem is reduced in this way 
to representation-directed
and to tame concealed algebras, but for both classes all the indecomposable modules are
accessible. This completes the proof of the Corollary. 
	    \bigskip
From now on, let $\Lambda$ be a non-distributive algebra and let $J$ be the radical of $\Lambda$.
Since the ideal lattice of $\Lambda$ is non-distributive, there are pairwise different 
ideals $I_0,\dots, I_3$ such
that $I_1\cap I_2 = I_2\cap I_3 = I_3 \cap I_1 = I_0$ and 
$I_1 + I_2 = I_2 + I_3 = I_3 + I_1.$ We can assume that $I_0 = 0,$ since with $\Lambda$ also
$\Lambda/I_0$ is non-distributive and the $\Lambda/I_0$-modules constructed can 
be considered as $\Lambda$-modules (annihilated by $I_0$).
Note that the existence of $I_3$ implies that the ideals $I_1$ and $I_2$ (considered as $\Lambda$-$\Lambda$-bimodules)
are isomorphic  and we can assume that these bimodules are simple bimodules. But since $\Lambda$
is a basic $k$-algebra and $k$ is algebraically closed, a simple $\Lambda$-$\Lambda$-bimodule $I$
is one-dimensional and there are primitive idempotents $e,f$ of $\Lambda$ 
(not necessarily different) such that $I = eIf.$
Thus, taking generators $\phi$ of $I_1$ and $\psi$ of $I_2$, these elements of $\Lambda$ are
linearly independent, there are primitive idempotents $e,f$ of $\Lambda$ such that $\phi = e\phi f$,
$\psi = e\psi f$ and $J\phi = J\psi = \phi J = \psi J = 0$ (conversely, the existence of such
elements $\phi, \psi\in \Lambda$ implies that $\Lambda$ is non-distributive).
	 
Let $E(e)$ be the injective envelope of the simple module $\Lambda e/J e$. In $E(e)$,
there are elements $x = fx,\  y = fy$ such that 
$$
 \phi x = 0,\quad u:= \psi x = \phi y \neq 0, \quad \psi y = 0.
$$
Note that $u$ is necessarily an element of the socle of $E(e)$.
Let  $V = \Lambda x+\Lambda y \subseteq E(e)$
     	
We consider direct sums of copies $V_{(i)} = V$, say $V^n = \bigoplus_{i=1}^n V_{(i)}$.
An element $v\in V$ will be denoted by $v_{(i)}$ when considered as an element of 
$V_{(i)} \subseteq V^n$. For $1\le i < n$
let $z_i = y_{(i)}+x_{(i+1)}.$
    
The following three submodules of $V^n$ (with $n\ge 1$) will be used:
$$
\align
  \Mminus &= \sum\nolimits_{i=1}^{n-1} \Lambda z_i,\quad \text{for}\quad n \ge 2,
 \quad \text{and} \quad  M(0)   = \Lambda u \subset V \cr
  R(n)   &= \Lambda x_{(1)}+ \Mminus, \cr
  W(n)   &= R(n)+ \Lambda y_{(n)}.
\endalign
$$
	\medskip
{\bf Proposition 1.} {\it The inclusions
$\Mminus \subset R(n)$ and $R(n) \subset W(n)$ are uniform.}
	 \medskip
The proof will use the following restriction lemma. Here, we denote by $B$ the subalgebra
of $\Lambda$ with basis $1, \phi,\psi$. It is a local algebra with radical square zero. 
If we consider a $\Lambda$-module $M$ as a $B$-module, then we write $_BM$.

   \medskip
{\bf Restriction Lemma 1.} {\it Let $M$ be a $\Lambda$-module. Assume that ${}_BM = N\oplus N'$ where $N$ is
an indecomposable  non-simple $B$-submodule and $N'$ is a semisimple $B$-module. Also, assume that 
$\soc {}_\Lambda M =\soc{}_B N$ (as vector spaces). Then $M$ is an indecomposable $\Lambda$-module.}
      \medskip
Proof. Let $M = M_1\oplus M_2$ be a direct decomposition of $M$ as a $\Lambda$-module, 
thus also ${}_BM = {}_B(M_1)\oplus {}_B(M_2)$. We apply
the theorem of Krull-Remak-Schmidt to the
direct decompositions $N\oplus N'= {}_BM = {}_B(M_1) \oplus {}_B(M_2)$ and see that
one of the summands $_B(M_1),{}_B(M_2)$, say $_B(M_1)$ can be written
in the form $N_1\oplus N'_1$ with $N_1$ isomorphic to $N$ and $N'_1$ semisimple
and then $_B(M_2)$ is also semisimple.
Since $N$ is an indecomposable non-simple $B$-module, we have $\soc N = \rad N.$ On the other hand, $\rad N' = 0 = \rad N'_1$ and also
$\rad {}_B(M_2) = 0.$ Thus,
$$
\align
 \soc M &= \soc N = \rad N = \rad N\oplus \rad N' =\rad {}_BM \cr 
 &= \rad N_1\oplus \rad N'_1\oplus \rad {}_B(M_2)
 = \rad N_1 \subseteq M_1.
\endalign
$$ 
But this implies that $M_2$ is zero (if $M_2 \neq 0$, then also $\soc M_2 \neq 0$
and of course $\soc M = \soc M_1\oplus \soc M_2).$
	\medskip
The indecomposable $B$-modules are well-known, since $B$ is stably 
equivalent to the Kronecker algebra $kQ$
(see for example [ARS], exercise X.3, or [Be], chapter 4.3; recall that 
the Kronecker quiver $Q$ is 
given by two vertices, say $a$ and $b$, and two arrows $a \to b$). 
For any $n>1$, there are up to isomorphism precisely indecomposable
$B$-modules of length $2n+1$, one is said to be {\it preprojective} (its socle
has length $n+1$, its top length $n$), the other one {\it preinjective}
(with socle of length $n$ and top of length $n+1$. The remaining non-simple 
indecomposables
are said to be {\it regular;} they have even length (and the length of the socle coincides
with the length of the top). For any $n\ge 1$, there is a up to isomorphism 
a unique indecomposable regular module of length $2n$ such that the kernel of
the multiplication by $\phi$ has dimension $n+1$.

    \medskip
Proof of proposition 1. 
We will consider $\Lambda$-modules $U$ with $\Mminus \subseteq U \subseteq V^n$; note that
for such a module $U$, one has $\soc U = \sum_{i=1}^{n} ku_{(i)}.$ 
Always, we will see that $_BU$ is the direct sum of
an indecomposable $B$-module $N$ and a semisimple $B$-module $N'.$
    \smallskip
{\bf(1)} {\it The inclusion $\Mminus \subseteq Jx_{(1)}+\Mminus$ is uniform for $n\ge 1.$}
	     
Proof. Consider a $\Lambda$-module $U$ with $\Mminus \subseteq U \subseteq Jx_{(1)}+\Mminus.$
If $n = 1$, then $U$ is a non-zero submodule of the uniform module $V$, thus indecomposable.
Let $n \ge 2.$ Let 
$$
 N = \sum_{i=1}^{n-1} Bz_i = \sum_{i=1}^{n-1} kz_i + \sum_{i=1}^{n} ku_{(i)},
$$ 
here we use that $\phi(z_{(i)}) = u_{(i)}$ and $\psi(z_{(i)}) = u_{(i+1)}$, for
$1\le i < n.$ Note that $N$ is the indecomposable preprojective $B$-module 
of length $2n-1 >1$ and its socle is $\soc {}_BN = 
\sum_{i=1}^{n} ku_{(i)}$.
Thus, we see that $\soc {}_BN = \soc U.$ 
On the other hand, $\Mminus = J\Mminus+N$, thus
$Jx_{(1)}+\Mminus = Jx_{(1)} + J\Mminus+N$. Since $\phi J= 0 = \psi J$, it follows that $Jx_{(1)} + J\Mminus$
is semisimple as a $B$-module. Thus $Jx_{(1)}+\Mminus$ is as a $B$-module the sum of $N$ and 
a semisimple $B$-module, and therefore also $U$ is as a $B$-module the sum of $N$ and 
a semisimple $B$-module $N'$. Altogether we see that we can apply the restriction lemma to the 
$\Lambda$-module $U$ and the $B$-modules $N,N'$ and conclude that $U$ is indecomposable.
	\smallskip	
{\bf(2)} {\it The inclusion $R(n) \subseteq Jy_{(n)}+R(n)$ is uniform for $n\ge 1.$}
	      
The proof is similar to that of (1), now
we consider a $\Lambda$-module $U$ with $R(n) \subseteq U \subseteq Jy_{(n)}+R(n)$ and can again
assume that $n \ge 2.$ This time, let 
$$
 N = Bx_{(1)}+\sum_{i=1}^{n-1} Bz_i = kx_{(1)}+ 
\sum_{i=1}^{n-1} kz_i + \sum_{i=1}^{n} ku_{(i)}.
$$ 
The $B$-module $N$ is regular indecomposable of length $2n >1$ and 
the kernel of the multiplication by $\phi$ has dimension $n+1$. The socle of $N$ is 
$\sum_{i=1}^{n} ku_{(i)} = \soc U$.
On the other hand, $R(n) = JR(n)+N$, thus
$Jy_{(n)}+R(n) = Jy_{(n)} + JR(n)+N$, and  $Jy_{(n)} + JR(n)$
is semisimple as a $B$-module. Since  $Jy_{(n)}+R(n)$ is as a $B$-module the sum of $N$ and 
a semisimple $B$-module, also ${}_BU$ is the sum of $N$ and 
a semisimple $B$-module $N'$. We apply again the restriction lemma to the 
$\Lambda$-module $U$ and the $B$-modules $N,N'$.
	\smallskip	
{\bf(3)} {\it The module $W(n)$ is indecomposable for $n\ge 1.$}
	      
The proof is again similar: let $U= W(n)$ and $n \ge 2.$ 
Now let 
$$
 N = Bx_{(1)}+ By_{(n)}+ \sum_{i=1}^{n-1} Bz_i = kx_{(1)}+ ky_{(n)}+ 
\sum_{i=1}^{n-1} kz_i + \sum_{i=1}^{n} ku_{(i)}.
$$
The $B$-module $N$ is the preinjective indecomposable $B$-module of length 
$2n+1 >1$, and its socle is $\sum_{i=1}^{n} ku_{(i)} = \soc U$.
On the other hand, $W(n) = JW(n)+N$, and $JW(n)$
is semisimple as a $B$-module. As before, we see that ${}_BU$ is  the sum of $N$ and 
a semisimple $B$-module $N'$. The restriction lemma shows that $U$ is indecomposable.
  \smallskip
{\bf (4)} {\it Let $M = M'+L$ be an indecomposable $\Lambda$-module with submodules $M'$ and $L$
such that $L$ is local. If $U$ is a $\Lambda$-module with $M'\subseteq U \subseteq M$, then
$U \subseteq M'+JL$ or else $U = M.$}
   
Proof. Let $U \subseteq M$ be a submodule
which is not contained in $M'+JL$. Then, in particular, 
$M'+JL$ is a proper submodule of $M$, and actually $M'+JL$ is a maximal submodule of $M$
(namely, the composition of the inclusion map $L \subseteq M = M'+L$ and the
projection $M \to M/(M'+JL)$ is surjective and contains $JL$ in its kernel, 
but $L/JL$ is simple). 

It follows that $M = U+(M'+JL) = U+JL$. 
Let $L = \Lambda m$ for some $m\in L$. Since $M = U+Jm,$
we see that $m = u+a m$ with $u\in U$ and $a\in J$, thus $(1-a)m = u\in U$. But
since $a\in J$, we know that $1-a$ is invertible in the ring $\Lambda$, 
therefore also $m\in U$.
As a consequence, $M = M'+\Lambda m \subseteq U$ and therefore $M = U.$
   \medskip
It follows from (2) that $R(n)$ is indecomposable, thus (1) and (4) show that the
inclusion $\Mminus\subset R(n)$ is uniform. Similarly, (2), (3) and (4) show that
the inclusion $R(n)\subset W(n)$ is uniform. This completes the proof of proposition 1.
    \bigskip
{\bf Proposition 2.} {\it For $n\ge 1$, there are couniform projections $M(n) \to R(n)$ and
$R(n+1) \to W(n)$.}
	\medskip
Proof: First, consider the embedding $R(n+1) \subset V^{n+1} = \bigoplus_{i=1}^{n+1}V_i$ and the
submodule $X = R(n+1)\cap V_{n+1} \subset R(n+1).$ Note that $R(n+1)/X = W(n)$, since for the canonical projection $R(n+1) \to R(n+1)/X$ we have $z_n \mapsto y_{(n)}$, whereas 
$x_{(1)} \mapsto x_{(1)},$ $z_i \mapsto z_i$ for $1 \le i \le n-1.$ 

Similarly, consider the embedding $M(n) \subset \bigoplus_{i=1}^{n+1}V_i$ and the
submodule $Y = M(n)\cap V_{(1)}\subset M(n).$ For the canonical projection
$M(n) \to M(n)/Y,$ we have $z_1 \mapsto x_{(2)}$, and 
$z_i \mapsto  z_{i+1}$ for $1 \le i \le n-1,$ thus we can identify
$M(n)/Y$ with $R(n)$ (where $R(n)$ is now considered as a submodule 
of $\bigoplus_{i=2}^{n+1} V_i$).
	 
In order to see that these projections $R(n+1) \to R(n+1)/X$ and $M(n) \to M(n)/Y$ are
couniform, we proceed as in the proof of Proposition 1, or better dually. In particular, 
we have to use the dual of
the restriction lemma 1 (here, 
instead of looking at the socles of ${}_\Lambda M$ and ${}_BN$, we assume that
the tops of ${}_\Lambda M$ and ${}_BN$ coincide):
   \medskip
{\bf Restriction Lemma 2.} {\it Let $M$ be a $\Lambda$-module. Assume that ${}_BM = N\oplus N'$ where $N$ is
an indecomposable  non-simple $B$-submodule and $N'$ is a semisimple $B$-module. Also, assume that
there is a vector subspace $T$ of $N$ such that  
$M = T\oplus \rad {}_\Lambda M$ and $N = T\oplus \rad {}_B N$ as vector spaces. 
Then $M$ is an indecomposable $\Lambda$-module.}
   \medskip
This completes the proof of proposition 2 and also that of the theorem.
     \medskip
{\bf Remark.} Note that in general the inclusion $\Mminus \subset W(n)$
is not uniform. Consider for $\Lambda$ the Kronecker algebra $kQ$, and look at
the submodules $U, U'$ of $W(2)$ generated by the elements 
$z = x_{(1)}+y_{(1)}+x_{(2)}+y_{(2)}$ and $z'= x_{(1)}-y_{(1)}-x_{(2)}+y_{(2)}$, respectively.
We have $\dim U = \dim U' = 2$. Assume now that the characteristic of $k$ is different from $2$.
Then $U \neq U'$ and even $U\cap U' = 0.$ Thus $U\oplus U'$ is a decomposable
submodule of $W(2)$. Also, $M(1)$ is contained in $U\oplus U'$ (as the submodule
generated by $\frac12(z-z')).$
		    \medskip
{\bf Acknowledgment.} The author is indebted to Dieter Vossieck for suggesting the
name {\it accessible}. 
     \medskip
{\bf References.}

\frenchspacing
\item{[ARS]} Auslander, M., Reiten, I., Smal\o{}, S.O.: Representation
 Theory of Artin Algebras. Cambridge University Press 1995.

\item{[Be]} Benson, D.J.: Representations and Cohomology I. Cambridge University
  Press 1991.

\item{[B1]} Bongartz, K.: Indecomposables over representation-finite algebras
  are extensions of an indecomposable and a simple. Math. Z. 187 (1984), 75-80.

\item{[B2]} Bongartz, K.: On degenerations and extensions of finite dimensional
 modules. Adv. Math. 121 (1996), 245-287.

\item{[B3]} Bongartz, K.: Indecomposables live in all smaller lengths. Preprint.
 arXiv:0904.4609

\item{[J]} Jans, J. P.: On the indecomposable representations of algebras.
Annals of Mathematics 66 (1957), 418-429.

\item{[R]}  A.V.Roiter, A.V.: Unboundedness of the dimension of the indecomposable representations
  of an algebra which has infinitely many indecomposable representations. 
  Izv. Akad. Nauk SSSR. Ser. Mat. 32 (1968), 1275-1282. 

	    \medskip

{\rmk Fakult\"at f\"ur Mathematik, Universit\"at Bielefeld \par
POBox 100\,131, \ D-33\,501 Bielefeld, Germany \par
e-mail: \ttk ringel\@math.uni-bielefeld.de \par}

\bye